# On Some Expansion Theorems Involving Confluent Hypergeometric ${}_2F_2(x)$ Polynomial

**Yashoverdhan Vyas*and Kalpana Fatawat**

Department of Mathematics, School of Engineering, Sir Padampat Singhania University Udaipur, Rajasthan.

**Abstract:** Recently, Rathie and Kılıçman (2014) employed Kummer-type transformation for ${}_2F_2(a, d+1; b, d\,;\, x)$ to develop certain classes of expansions theorems for ${}_2F_2(x)$ hypergeometric polynomial. Our aim is to deduce Kummer-type transformation for ${}_2F_2(a, d+2; b, d\,;\, x)$ and utilize it to develop some new expansion theorems for the confluent hypergeometric ${}_2F_2(x)$ polynomial. We also obtain a well-known result given by Kim et al. (Integral Transforms Spec. Funct. 23(6); 435-444, 2012) and many other new results as particular cases of our theorems.

## INTRODUCTION

The subject “Special Functions” is an important branch of “Analysis” and most of the “Special Functions” are solutions of “Differential Equations”. A number of problems from different branches of science and engineering can be expressed in the form of differential equations satisfied by the “Special Functions”. Further, there is no other way better than the hypergeometric type functions to generalize and unify all the “Special Functions”, orthogonal polynomials and other well-known elementary functions.

The generalized hypergeometric function with arbitrary number of numerator and denominator parameters is defined as follows:

$$ {}_pF_q\left[\begin{matrix} a_1, a_2, \ldots\ldots a_p \\ b_1, b_2, \ldots\ldots b_q \end{matrix}; z\right] = \sum_{n=0}^{\infty} \frac{(a_1)_n (a_2)_n \ldots\ldots (a_p)_n}{(b_1)_n (b_2)_n \ldots\ldots (b_p)_n} \frac{z^n}{n!} \tag{1.1} $$

It converges for every finite value of $z$ if $p \leq q$, diverges for every $z, z \neq 0$, if $p > q+1$, converges for all $|z|<1$, if $p = q+1$ and also converges for $|z| = 1$ with $p = q+1$, if $\mathrm{Re}(\sum_{j=0}^{q} b_j - \sum_{j=0}^{p} a_j) > 0$, provided $b_j$ is neither zero and nor a negative integer.

The ${}_pF_q(z) = w$ defined in equation (1.1) satisfies following differential equation:

$$ \left[\theta \prod_{j=1}^{q}(\theta + b_j - 1) - z\prod_{i=1}^{p}(\theta + a_i)\right] w = 0, \text{ where } \theta \equiv z\frac{d}{dz} \tag{1.2} $$

In this paper, we will use confluent hypergeometric series ${}_pF_p$ and the generalized hypergeometric series ${}_pF_q$ with $p < q$, which are always convergent and the ${}_{q+1}F_q(\pm 1)$ type of series for which the convergence condition is mentioned above. For additional information on the importance, development and applications of summation, transformation, expansion and other types of theorems on hypergeometric type functions, we refer the books by Rainville [1], Slater[2], Prudnikov *et al.*[3], Brychkov[4] and Olver *et al.*[5].

Further, in many branches of pure and applied mathematics, the Laguerre polynomial, which is a terminating form of the confluent hypergeometric function ${}_1F_1$ defined by

$$ L_n^{(v)}(x) = \frac{(v+1)_n}{n!} {}_1F_1\left[\begin{matrix} -n \\ v+1 \end{matrix}; x\right] \tag{1.3} $$

occurs frequently; see, for example, Erdelyi[6].

Kim *et al.*[7] established a general transformation involving the generalized hypergeometric function by the method of series manipulation. The well-known Kummer's first theorem, the classical Gauss summation theorem and the generalized Kummer summation theorem were then applied to obtain a new class of summation or expansion formulae involving the Laguerre polynomial, which did not appear previously in the literature.

Recently, Brychkov[8] and Kim *et al.*[9] have further studied and developed additional results on expansion involving Laguerre polynomials, which clearly indicate the current importance and interest about these types of results.

Rathie *et al.*[10] moved one step ahead and developed such expansions for ${}_2F_2$ [a, d+1; b, d; x] by using extensions of some contiguous Kummer and Gauss summation theorems and a Kummer type transformation discovered by Paris[11].

$$e^{-x}\,{}_2F_2\left[\begin{matrix} a, d+1 \\ b+1,\ d \end{matrix};\ x\right] = {}_2F_2\left[\begin{matrix} b-a,\ f+1 \\ b+1,\ f \end{matrix};\ -x\right] \text{ Where, } f = \frac{d(a-b)}{(a-d)} \qquad (1.4)$$

Rathie and Paris[11] have recently, developed expansions involving Laguerre polynomials with some additional coefficient terms. But, the expansions for a general $_2F_2[x]$ have not been studied, probably, because of the difficulty mentioned in Paris[11] on developing Kummer type transformations for a general $_2F_2[x]$.

In this paper, it is shown, how another transformation for $_2F_2$ [a, d+2; b, d; x] can be developed and utilized to find certain additional expansions involving a $_2F_2$ confluent hypergeometric polynomial. On particularization, our results not only lead to a previously known result by Kim *et al.*[7] but also produce expansions involving two Laguerre polynomials.

In the next three sections the statements and derivations of some required extensions of contiguous Kummer summation theorems and Gauss summation theorem and a Kummer-type transformation for confluent hypergeometric function $_2F_2(x)$ and a general transformation formula to derive the expansion theorems are given. Then, in the last two sections, the main expansion theorems and their special cases are discussed.

## SOME EXTENSIONS OF CONTIGUOUS KUMMER AND GAUSS SUMMATION THEOREMS

$${}_3F_2\left[\begin{matrix} a,b,d+2 \\ a-b,\ d \end{matrix};\ -1\right] = \frac{\Gamma(\frac{1}{2})\Gamma(a-b)}{2^a}\left\{\frac{(1-\frac{2b}{d}+\frac{b(a+1)}{2d(d+1)})}{\Gamma(\frac{a}{2})\Gamma(\frac{a}{2}-b+\frac{1}{2})} + \frac{(1-\frac{ab}{2d(d+1)})}{\Gamma(\frac{a}{2}+\frac{1}{2})\Gamma(\frac{a}{2}-b)}\right\} \qquad (2.1)$$

$${}_3F_2\left[\begin{matrix} a,b,d+2 \\ 3+a-b,\ d \end{matrix};\ -1\right] = \frac{\Gamma(\frac{1}{2})\Gamma(3+a-b)}{2^a(1-b)(2-b)}\left\{\frac{\left\{-2+\frac{2(2a-b+2)}{d}-\frac{2(a-b+2)(a+1)}{d(d+1)}\right\}}{\Gamma(\frac{a}{2})\Gamma(\frac{a}{2}-b+\frac{3}{2})} + \frac{\left\{(a-b+1)+\frac{a(-2a+3b-4)}{d}+\frac{a\{2(a-b+3)(a-b+1)-(b+1)(b-2)\}}{2d(d+1)}\right\}}{\Gamma(\frac{a}{2}+\frac{1}{2})\Gamma(\frac{a}{2}-b+2)}\right\} \qquad (2.2)$$

$${}_3F_2\left[\begin{matrix} a,b,d+2 \\ c+1,\ d \end{matrix};\ 1\right] = \frac{\Gamma(c+1)\Gamma(c-a-b)}{\Gamma(c-a+1)\Gamma(c-b+1)}\left\{(c-a-b)+\frac{2ab}{d}+\frac{ab(a+1)(b+1)}{d(d+1)(c-a-b-1)}\right\} \qquad (2.3)$$

**Proof:**

Denoting the left-hand side of (2.1) by S and express $_3F_2(-1)$ as a series, after some simplification, we get

$$S=\sum_{n\geq 0}\frac{(a)_n(b)_n}{(a-b)_n}\left(1+\frac{n}{d}\right)\left(1+\frac{n}{(d+1)}\right)\frac{(-1)^n}{n!}$$

After, simplifying the above expression and using contiguous extension of Kummer transformation and other formulae listed in Rathie *et al.*[10], followed by some more simplifications, we get the above result (2.1). Likewise, two more contiguous relations given by the equations $(2.2)$ and $(2.3)$ can be derived.

## A KUMMER-TYPE TRANSFORMATION FOR ${}_2F_2\,(a,d+2;b,d;\ x)$

$${}_2\mathrm{F}_2\left[\begin{matrix}a,d+2\\ b,\ d\end{matrix};\ x\right]=\sum_{n\geq 0}\frac{(a)_n}{(b)_n}\left(1+\frac{2n}{d}+\frac{n(n-1)}{d(d+1)}\right)\frac{x^n}{n!}$$

Simplifying right hand side and using equation (5) of Paris[11], we obtain

$${}_2\mathrm{F}_2\left[\begin{matrix}a,d+2\\ b,\ d\end{matrix};\ x\right]=\left\{{}_2F_2\left[\begin{matrix}a,d+1\\ b,d\end{matrix};\ x\right]+\frac{ax}{bd}\,{}_2F_2\left[\begin{matrix}a+1,d+2\\ b+1,d+1\end{matrix};\ x\right]\right\}$$

Now, using equation (1.4), we get the required extension of Kummer-type transformation in the following form :

$${}_2\mathrm{F}_2\left[\begin{matrix}a,d+2\\ b,\ d\end{matrix};\ x\right]=e^x\left\{{}_2F_2\left[\begin{matrix}b-a,\ f+1\\ b+1,\ f\end{matrix};\ -x\right]\right.$$
$$\left.+\frac{ax}{(b+1)d}\,{}_2F_2\left[\begin{matrix}b-a,\ f'+1\\ b+2,\ f'\end{matrix};\ -x\right]\right\}\qquad(3.1)$$

where,

$$f=\frac{d(a-b)}{(a-d)},\qquad f'=\frac{(d+1)(a-b)}{(a-d)}\qquad(3.2)$$

## A GENERAL TRANSFORMATION FORMULA

The needed transformation formula involving the confluent hypergeometric function and the generalized hypergeometric function is :

$$\sum_{n\geq 0}\frac{(a_1)_n......(a_p)_n}{(b_1)_n......(b_q)_n}\frac{(d)_n x^n y^n}{n!}\,{}_2F_2\left[\begin{matrix} d+n, e+2 \\ f+1, e\end{matrix}; x\right]=$$
$$\sum_{n\geq 0}\frac{(d)_n (e+2)_n}{(f+1)_n (e)_n}\frac{x^n}{n!}\,{}_{p+3}F_{q+1}\left[\begin{matrix} -n, 1-e-n, -f-n, a_1......a_p \\ -1-e-n, b_1......b_q\end{matrix}; y\right] \qquad (3.3)$$

**Proof**: Denote left-hand side of $(3.3)$ by S and express ${}_2F_2(x)$ as a series to obtain

$$S=\sum_{n\geq 0}\frac{(a_1)_n......(a_p)_n}{(b_1)_n......(b_q)_n}\frac{(d)_n x^n y^n}{n!}\sum_{m\geq 0}\frac{(d+n)_m (e+2)_m}{(f+1)_m (e)_m}\frac{x^m}{m!}$$

By applying simplifications of shifted factorials and series manipulation technique

$$\sum_{k,n\geq 0} A(k,n)=\sum_{n\geq 0}\sum_{k=0}^{n} A(k,n-k)$$

, we get equation (3.3).

## EXPANSION THEOREMS INVOLVING THE CONFLUENT HYPERGEOMETRIC ${}_2F_2$ POLYNOMIAL

### Theorem 1

$$e^{-x}\sum_{n\geq 0}\frac{x^n}{n!}\left\{{}_2F_2\left[\begin{matrix} -n, g+1 \\ v+2, g\end{matrix}; x\right]+\frac{(v+1+n)(-x)}{e(v+2)}\,{}_2F_2\left[\begin{matrix} -n, g'+1 \\ v+3, g'\end{matrix}; x\right]\right\}=$$
$${}_3F_4\left[\begin{matrix} \frac{v}{2}+\frac{1}{2}, \frac{e}{2}+1, \frac{e}{2}+\frac{3}{2} \\ v+1, \frac{v}{2}+\frac{3}{2}, \frac{e}{2}+\frac{1}{2}, \frac{e}{2}\end{matrix}; -x^2\right]+\frac{x^2}{e(e+1)}\,{}_0F_1\left[\begin{matrix} - \\ v+2\end{matrix}; -x^2\right]+$$
$$\frac{(e+2)x}{(v+2)e}\,{}_3F_4\left[\begin{matrix} \frac{v}{2}+1, \frac{e}{2}+2, \frac{e}{2}+\frac{3}{2} \\ v+2, \frac{v}{2}+2, \frac{e}{2}+\frac{1}{2}, \frac{e}{2}+1\end{matrix}; -x^2\right]-\frac{2x}{e}\,{}_1F_2\left[\begin{matrix} \frac{e}{2}+\frac{3}{2} \\ v+2, \frac{e}{2}+\frac{1}{2}\end{matrix}; -x^2\right]$$
$$-\frac{x^3}{e(e+1)(v+2)}\,{}_0F_1\left[\begin{matrix} - \\ v+3\end{matrix}; -x^2\right] \qquad (4.1)$$

Where,

$$g=\frac{en}{v+1+n-e} \qquad \text{And} \qquad g'=\frac{(e+1)n}{v+1+n-e}$$

**Proof**: Applying the extension of Kummer-type transformation $(3.1)$ to left hand side of $(3.3)$ we obtain

$$e^x \sum_{n\geq 0} \frac{(a_1)_n ......(a_q)_n (d)_n}{(b_1)_n ......(b_q)_n} \frac{x^n y^n}{n!} \left\{ {}_2F_2 \left[ \begin{matrix} f-d-n, g+1 \\ f+1,\ g \end{matrix} ; -x \right] \right.$$

$$\left. + \frac{(d+n)x}{e(f+1)} \ {}_2F_2 \left[ \begin{matrix} f-d-n, g'+1 \\ f+2,\ g' \end{matrix} ; -x \right] \right\}$$

$$= \sum_{n\geq 0} \frac{(d)_n (e+2)_n}{(f+1)_n (e)_n} \frac{x^n}{n!} \ {}_{p+3}F_{q+1} \left[ \begin{matrix} -n, 1-e-n, -f-n, a_1......a_p \\ -1-e-n,\ b_1......b_q \end{matrix} ; y \right]$$

(4.2)

where $g = \dfrac{e(d+n-f)}{d+n-e}$ and $g' = \dfrac{(e+1)(d+n-f)}{d+n-e}$ .

Substituting $d = f$ in (4.2), we find

$$e^x \sum_{n\geq 0} \frac{(a_1)_n ......(a_q)_n (f)_n}{(b_1)_n ......(b_q)_n} \frac{x^n y^n}{n!} \left\{ {}_2F_2 \left[ \begin{matrix} -n, g+1 \\ f+1,\ g \end{matrix} ; -x \right] + \frac{(f+n)x}{e(f+1)} \ {}_2F_2 \left[ \begin{matrix} -n, g'+1 \\ f+2,\ g' \end{matrix} ; -x \right] \right\} =$$

$$\sum_{n\geq 0} \frac{(f)_n (e+2)_n}{(f+1)_n (e)_n} \frac{x^n}{n!} \ {}_{p+3}F_{q+1} \left[ \begin{matrix} -n, 1-e-n, -f-n, a_1......a_p \\ -1-e-n,\ b_1......b_q \end{matrix} ; y \right]$$

and replacing $x$ with $-x$ and setting $f = v+1$, we obtain

$$e^{-x} \sum_{n\geq 0} \frac{(a_1)_n ......(a_q)_n (v+1)_n}{(b_1)_n ......(b_q)_n} \frac{(-x)^n y^n}{n!} \left\{ {}_2F_2 \left[ \begin{matrix} -n, g+1 \\ v+2,\ g \end{matrix} ; x \right] - \right.$$

$$\left. \frac{(v+1+n)x}{e(f+1)} \ {}_2F_2 \left[ \begin{matrix} -n, g'+1 \\ v+3,\ g' \end{matrix} ; x \right] \right\} =$$

$$\sum_{n\geq 0} \frac{(v+1)_n (e+2)_n}{(v+2)_n (e)_n} \frac{(-x)^n}{n!} \ {}_{p+3}F_{q+1} \left[ \begin{matrix} -n, 1-e-n, -v-1-n, a_1......a_p \\ -1-e-n,\ b_1......b_q \end{matrix} ; y \right]$$

(4.3)

Selecting $p = 0, q = 1, b_1 = v+1$ and $y = -1$ in (4.3), we arrive at

$$e^{-x} \sum_{n\geq 0} \frac{x^n}{n!} \left\{ {}_2F_2 \left[ \begin{matrix} -n, g+1 \\ v+2,\ g \end{matrix} ; x \right] - \frac{(v+1+n)x}{e(f+1)} \ {}_2F_2 \left[ \begin{matrix} -n, g'+1 \\ v+3,\ g' \end{matrix} ; x \right] \right\} =$$

$$\sum_{n\geq 0} \frac{(v+1)_n (e+2)_n}{(v+2)_n (e)_n} \frac{(-x)^n}{n!} \ {}_3F_2 \left[ \begin{matrix} -n, 1-e-n, -v-1-n \\ -1-e-n,\ v+1 \end{matrix} ; -1 \right]$$

(4.4)

Applying (2.1) to solve ${}_3F_2(-1)$ on the right hand side of (4.4), the resultant expression becomes

$$= \sum_{n\ge 0} \frac{(e+2)_n (v+1)_n}{(v+2)_n (e)_n} \frac{(-x)^n}{n!} \frac{\Gamma(\frac{1}{2})\Gamma(v+1)}{2^{-n}} \left\{ \frac{(1-\frac{2(v+1+n)}{1+e+n}+\frac{(v+1+n)(n-1)}{2(1+e+n)(e+n)})}{\Gamma(\frac{-n}{2})\Gamma(\frac{n}{2}+v+\frac{3}{2})} + \frac{(1-\frac{n(v+1+n)}{2(1+e+n)(e+n)})}{\Gamma(\frac{-n}{2}+\frac{1}{2})\Gamma(\frac{n}{2}+v+1)} \right\} \tag{4.5}$$

Now, we separate the terms of equation $(4.5)$ into even and odd powers of $x$. The desired result is obtained after some simplifications using gamma functions and shifted factorial identities.

**Theorem 2**

$$e^{-x} \sum_{n\ge 0} \frac{(v+1)_n}{(2-v)_n} \frac{x^n}{n!} \left\{ {}_2F_2\left[\begin{matrix} -n, g+1 \\ v+2,\ g \end{matrix}; x\right] - \frac{(v+1+n)x}{e(v+2)} {}_2F_2\left[\begin{matrix} -n, g'+1 \\ v+3,\ g' \end{matrix}; x\right] \right\} =$$

$$\frac{\Gamma(\frac{1}{2})\Gamma(2-v)}{2^{-v}\Gamma(-\frac{v}{2}+\frac{3}{2})\Gamma(-\frac{v}{2})} \left\{ -v\ {}_5F_6\left[\begin{matrix} 1, \frac{v}{2}+1, \frac{v}{2}+\frac{1}{2}, \frac{e}{2}+1, \frac{e}{2}+\frac{3}{2} \\ 2, \frac{3}{2}, \frac{v}{2}+\frac{3}{2}, \frac{e}{2}, \frac{e}{2}+\frac{1}{2}, -\frac{v}{2}+\frac{3}{2} \end{matrix}; -x^2\right] + \right.$$

$$\frac{2(v^2-1)}{(1+e)}\ {}_4F_5\left[\begin{matrix} 1, \frac{v}{2}+1, \frac{e}{2}+1, 2-v \\ 2, \frac{3}{2}, \frac{e}{2}, 1-v, -\frac{v}{2}+\frac{3}{2} \end{matrix}; -x^2\right] - \frac{v(v+1)(v-2)}{e(1+e)}\ {}_2F_3\left[\begin{matrix} 1, \frac{v}{2}+1 \\ 2, \frac{3}{2}, -\frac{v}{2}+\frac{3}{2} \end{matrix}; -x^2\right] -$$

$$\left. \frac{1}{e(1+e)}\ {}_3F_4\left[\begin{matrix} \frac{1}{2}, \frac{v}{2}+1, \frac{v}{2}+\frac{1}{2} \\ \frac{3}{2}, \frac{-1}{2}, \frac{v}{2}+\frac{3}{2}, -\frac{v}{2}+\frac{3}{2} \end{matrix}; -x^2\right] \right\} + \frac{\Gamma(\frac{1}{2})\Gamma(2-v)}{2^{-v-1}\Gamma(-\frac{v}{2}-\frac{1}{2})\Gamma(-\frac{v}{2}+1)}$$

$$\left\{ -\ {}_4F_5\left[\begin{matrix} 1, \frac{v}{2}+\frac{1}{2}, \frac{e}{2}+1, \frac{e}{2}+\frac{3}{2} \\ 2, \frac{3}{2}, \frac{e}{2}, \frac{e}{2}+\frac{1}{2}, -\frac{v}{2}+1 \end{matrix}; -x^2\right] + \frac{2v}{(1+e)}\ {}_4F_5\left[\begin{matrix} 1, \frac{v}{2}+\frac{1}{2}, \frac{e}{2}+1, v+1 \\ 2, \frac{3}{2}, \frac{e}{2}, v, -\frac{v}{2}+1 \end{matrix}; -x^2\right] - \right.$$

$$\left. \frac{v(v-1)}{e(1+e)}\ {}_3F_4\left[\begin{matrix} 1, \frac{v}{2}+\frac{1}{2}, \frac{v}{2}+1 \\ 2, \frac{3}{2}, \frac{v}{2}, -\frac{v}{2}+1 \end{matrix}; -x^2\right] \right\} + \frac{x\Gamma(\frac{1}{2})\Gamma(2-v)}{2^{-v}\Gamma(-\frac{v}{2}+2)\Gamma(-\frac{v}{2}-\frac{1}{2})}$$

$$\left\{ \frac{2v(v+1)(e+2)}{3e(v+2)}\ {}_5F_6\left[\begin{matrix} 1, \frac{v}{2}+1, \frac{v}{2}+\frac{3}{2}, \frac{e}{2}+2, \frac{e}{2}+\frac{3}{2} \\ 2, \frac{5}{2}, \frac{v}{2}+2, \frac{e}{2}+1, \frac{e}{2}+\frac{1}{2}, -\frac{v}{2}+2 \end{matrix}; -x^2\right] - \right.$$

$$\frac{2(2v-3)(v+1)}{3e}\ {}_4F_5\left[\begin{matrix} 1, \frac{v}{2}+\frac{3}{2}, \frac{e}{2}+\frac{3}{2}, -v+\frac{5}{2} \\ 2, \frac{5}{2}, \frac{e}{2}+\frac{1}{2}, -v+\frac{3}{2}, -\frac{v}{2}+2 \end{matrix}; -x^2\right] + \frac{2v(v+1)(v-2)}{3e(1+e)}\ {}_2F_3\left[\begin{matrix} 1, \frac{v}{2}+\frac{3}{2} \\ 2, \frac{5}{2}, -\frac{v}{2}+2 \end{matrix}; -x^2\right] -$$

$$\left. \frac{2(v+3)(v+1)}{3e(1+e)(4-v^2)}\ {}_3F_4\left[\begin{matrix} 2, \frac{v}{2}+2, \frac{v}{2}+\frac{5}{2} \\ 3, \frac{5}{2}, \frac{v}{2}+3, -\frac{v}{2}+3 \end{matrix}; -x^2\right] \right\} + \frac{x\Gamma(\frac{1}{2})\Gamma(2-v)}{2^{-v-1}\Gamma(-\frac{v}{2}-1)\Gamma(-\frac{v}{2}+\frac{3}{2})}$$

$$\left\{ \frac{2(v+1)(e+2)}{3e(v+2)}\ {}_4F_5\left[\begin{matrix} 1, \frac{v}{2}+1, \frac{e}{2}+2, \frac{e}{2}+\frac{3}{2} \\ 2, \frac{5}{2}, \frac{e}{2}+1, \frac{e}{2}+\frac{1}{2}, -\frac{v}{2}+\frac{3}{2} \end{matrix}; -x^2\right] - \right.$$

$$\left. \frac{2(2v+1)(v+1)}{3e(v+2)}\ {}_4F_5\left[\begin{matrix} 1, \frac{v}{2}+1, \frac{e}{2}+\frac{3}{2}, v+\frac{3}{2} \\ 2, \frac{5}{2}, \frac{e}{2}+\frac{1}{2}, v+\frac{1}{2}, -\frac{v}{2}+\frac{3}{2} \end{matrix}; -x^2\right] + \frac{2(v+1)^2(v-1)}{3e(1+e)(v+2)}\ {}_3F_4\left[\begin{matrix} 1, \frac{v}{2}+1, \frac{v}{2}+\frac{3}{2} \\ 2, \frac{5}{2}, -\frac{v}{2}+\frac{3}{2}, \frac{v}{2}+\frac{1}{2} \end{matrix}; -x^2\right] \right\}\Bigg\}$$

Where $g=\frac{en}{v+1+n-e}$ and $g'=\frac{(e+1)n}{v+1+n-e}$ (4.6)

**Proof:** Selecting $p=0, q=1, b_1=2-v$ and $y=-1$ in $(4.3)$, we get

$$e^{-x}\sum_{n\geq 0}\frac{(1-v)_n\, x^n}{(2-v)_n\, n!}\left\{{}_2F_2\left[\begin{matrix}-n, g+1\\ v+2,\ g\end{matrix};x\right]-\frac{(v+1+n)x}{e(v+2)}\,{}_2F_2\left[\begin{matrix}-n, g'+1\\ v+3,\ g'\end{matrix};x\right]\right\}=$$
$$\sum_{n\geq 0}\frac{(v+1)_n(e+2)_n}{(v+2)_n(e)_n}\frac{(-x)^n}{n!}\,{}_3F_2\left[\begin{matrix}-n, 1-e-n,\ -v-1-n\\ -1-e-n,\ 2-v\end{matrix};-1\right] \quad (4.7)$$

Next, we apply equation $(2.2)$ to solve ${}_3F_2(-1)$ on the right hand side of $(4.7)$, to have

$$\text{RHS of }(4.7)=\sum_{n\geq 0}\frac{(e+2)_n(v+1)_n}{(v+2)_n(e)_n}\frac{(-x)^n}{n!}\frac{\Gamma(\frac{1}{2})\Gamma(2-v)}{2^{-v-1-n}(2+n)(n+1)}$$

$$\left\{\frac{\left(-v+\frac{(2v-2-n)(v+1+n)}{1+e+n}-\frac{(v-2)v(v+1+n)}{(1+e+n)(e+n)}+\frac{(n-1)(2+n)}{2(1+e+n)(e+n)}\right)}{\Gamma(\frac{-n}{2}-\frac{v}{2})\Gamma(\frac{n}{2}-\frac{v}{2}+\frac{3}{2})}+\frac{2\left(-1+\frac{(2v+n)}{(1+e+n)}-\frac{(v-1)(v+n)}{2(1+e+n)(e+n)}\right)}{\Gamma(\frac{-n}{2}-\frac{v}{2}-\frac{1}{2})\Gamma(\frac{n}{2}-\frac{v}{2}+1)}\right\} \quad (4.8)$$

Now, separating the terms appearing on the right hand side of $(4.8)$ into even and odd power of x. After some simplifications using gamma functions and shifted factorial identities and following the steps similar to those mentioned in theorem 1, we obtain theorem 2.

**Theorem 3:**

$$e^{-x}\sum_{n\geq 0}\frac{(v+1)_n(-x)^n}{(\mu)_n\, n!}\left\{{}_2F_2\left[\begin{matrix}-n, g+1\\ v+2,\ g\end{matrix};x\right]+\frac{(v+1+n)(-x)}{e(v+2)}\,{}_2F_2\left[\begin{matrix}-n, g'+1\\ v+3,\ g'\end{matrix};x\right]\right\}=$$
$${}_4F_4\left[\begin{matrix}v+1, e+2, \frac{\mu}{2}+\frac{v}{2}+1, \frac{\mu}{2}+\frac{v}{2}+\frac{1}{2}\\ v+2, e, \mu, \mu+v+1\end{matrix};-4x\right]+\frac{2x(v+1)}{e\mu}\,{}_3F_3\left[\begin{matrix}e+2, \frac{\mu}{2}+\frac{v}{2}+1, \frac{\mu}{2}+\frac{v}{2}+\frac{3}{2}\\ e+1, \mu+1, \mu+v+2\end{matrix};-4x\right]+$$
$$\frac{(v+1)(v+2)x^2}{e(e+1)\mu(\mu+1)}\,{}_3F_3\left[\begin{matrix}v+3, \frac{\mu}{2}+\frac{v}{2}+2, \frac{\mu}{2}+\frac{v}{2}+\frac{3}{2}\\ v+2, \mu+2, \mu+v+3\end{matrix};-4x\right]$$

Where,

$g=\frac{en}{v+1+n-e}$ and

$g'=\frac{(e+1)n}{v+1+n-e}$ (4.9)

**Proof:** To derive $(4.9)$, let $p=0, q=1, b_1=\mu$ and $y=1$ in equation $(4.3)$, we use the Kummer extension given in equation (2.3) and proceeding on the line of the proofs above two theorems, we obtain the required result.

**SPECIAL CASES:**

For $e = v+1$, in (4.9), we get $g = v+1$ and $g' = v+2$.

And so for $e = v+1$, Theorem 1 yields

$$e^{-x}\sum_{n\geq 0}\frac{x^n}{(1+v)_n}L_n^{(v+1)}(x) = \Gamma(v+1)x^{-v-1}\{(x+v+1)J_{v+1}(2x) - xJ_{v+2}(2x)\} \tag{5.1}$$

Replacing $v+1$ by $v$ in $(5.1)$, we obtain a well-known result investigated by Kim *et al.*[7]

$$e^{-x}\sum_{n\geq 0}\frac{x^n}{(v)_n}L_n^{(v)}(x) = \Gamma(v)x^{1-v}\{J_{v-1}(2x) + xJ_v(2x)\} \tag{5.2}$$

By choosing $v = \pm ½$ in above result, special summation is obtained in terms of trigonometric functions. For further information, we refer to Kim *et al.*[7]

Similarly, for $e = v+1$, Theorem 2 yields,

$$e^{-x}\sum_{n\geq 0}\frac{x^n}{(2-v)_n}\left\{L_n^{(v)}(x) - \tfrac{x}{v+2}L_n^{(v+1)}(x)\right\}$$

$$= \frac{\Gamma(\tfrac{1}{2})\Gamma(2-v)}{2^{-v}\Gamma(-\tfrac{v}{2}+\tfrac{3}{2})\Gamma(-\tfrac{v}{2})}\left\{-v\ {}_2F_3\left[\begin{matrix}1, \tfrac{v}{2}+2 \\ 2, \tfrac{3}{2}, -\tfrac{v}{2}+\tfrac{3}{2}\end{matrix}; -x^2\right] + \tfrac{2(v^2-1)}{(2+v)}\ {}_4F_5\left[\begin{matrix}1, \tfrac{v}{2}+1, \tfrac{v}{2}+\tfrac{3}{2}, 2-v \\ 2, \tfrac{3}{2}, \tfrac{v}{2}+\tfrac{1}{2}, 1-v, -\tfrac{v}{2}+\tfrac{3}{2}\end{matrix}; -x^2\right]\right.$$

$$\left. - \tfrac{v(v-2)}{(v+2)}\ {}_2F_3\left[\begin{matrix}1, \tfrac{v}{2}+1 \\ 2, \tfrac{3}{2}, -\tfrac{v}{2}+\tfrac{3}{2}\end{matrix}; -x^2\right] - \tfrac{1}{(v+1)(2+v)}\ {}_3F_4\left[\begin{matrix}\tfrac{1}{2}, \tfrac{v}{2}+1, \tfrac{v}{2}+\tfrac{1}{2} \\ \tfrac{3}{2}, \tfrac{-1}{2}, \tfrac{v}{2}+\tfrac{3}{2}, -\tfrac{v}{2}+\tfrac{3}{2}\end{matrix}; -x^2\right]\right\}$$

$$+ \frac{\Gamma(\tfrac{1}{2})\Gamma(2-v)}{2^{-v-1}\Gamma(-\tfrac{v}{2}-\tfrac{1}{2})\Gamma(-\tfrac{v}{2}+1)}\left\{-\ {}_3F_4\left[\begin{matrix}1, \tfrac{v}{2}+\tfrac{3}{2}, \tfrac{v}{2}+2 \\ 2, \tfrac{3}{2}, \tfrac{v}{2}+1, -\tfrac{v}{2}+1\end{matrix}; -x^2\right] + \tfrac{2v}{(v+2)}\ {}_3F_4\left[\begin{matrix}1, \tfrac{v}{2}+\tfrac{3}{2}, v+1 \\ 2, \tfrac{3}{2}, v, -\tfrac{v}{2}+1\end{matrix}; -x^2\right] -\right.$$

$$\left. \tfrac{v(v-1)}{(v+1)(v+2)}\ {}_3F_4\left[\begin{matrix}1, \tfrac{v}{2}+\tfrac{1}{2}, \tfrac{v}{2}+1 \\ 2, \tfrac{3}{2}, \tfrac{v}{2}, -\tfrac{v}{2}+1\end{matrix}; -x^2\right]\right\} + \frac{x\Gamma(\tfrac{1}{2})\Gamma(2-v)}{2^{-v}\Gamma(-\tfrac{v}{2}+2)\Gamma(-\tfrac{v}{2}-\tfrac{1}{2})}$$

$$\left\{\tfrac{2v(v+3)}{3(v+2)}\ {}_2F_3\left[\begin{matrix}1, \tfrac{v}{2}+\tfrac{5}{2} \\ 2, \tfrac{5}{2}, -\tfrac{v}{2}+2\end{matrix}; -x^2\right] - \tfrac{2(2v-3)}{3}\ {}_4F_5\left[\begin{matrix}1, \tfrac{v}{2}+\tfrac{3}{2}, \tfrac{v}{2}+2, -v+\tfrac{5}{2} \\ 2, \tfrac{5}{2}, \tfrac{v}{2}+1, -v+\tfrac{3}{2}, -\tfrac{v}{2}+2\end{matrix}; -x^2\right]\right.$$

$$\left. + \tfrac{2v(v-2)}{3(v+2)}\ {}_2F_3\left[\begin{matrix}1, \tfrac{v}{2}+\tfrac{3}{2} \\ 2, \tfrac{5}{2}, -\tfrac{v}{2}+2\end{matrix}; -x^2\right] - \tfrac{2(v+3)}{3(v+2)(4-v^2)}\ {}_3F_4\left[\begin{matrix}2, \tfrac{v}{2}+2, \tfrac{v}{2}+\tfrac{5}{2} \\ 3, \tfrac{5}{2}, \tfrac{v}{2}+3, -\tfrac{v}{2}+3\end{matrix}; -x^2\right]\right\}$$

$$+ \frac{x\Gamma(\tfrac{1}{2})\Gamma(2-v)}{2^{-v-1}\Gamma(-\tfrac{v}{2}-1)\Gamma(-\tfrac{v}{2}+\tfrac{3}{2})}\left\{\tfrac{2((v+3)}{3(v+2)}\ {}_4F_5\left[\begin{matrix}1, \tfrac{v}{2}+\tfrac{5}{2} \\ 2, \tfrac{5}{2}, -\tfrac{v}{2}+\tfrac{3}{2}\end{matrix}; -x^2\right] -\right.$$

$$\left. \tfrac{2(2v+1)}{3(v+2)}\ {}_3F_4\left[\begin{matrix}1, \tfrac{v}{2}+2, v+\tfrac{3}{2} \\ 2, \tfrac{5}{2}, v+\tfrac{1}{2}, -\tfrac{v}{2}+\tfrac{3}{2}\end{matrix}; -x^2\right] + \tfrac{2(v+1)(v-1)}{3(v+2)^2}\ {}_3F_4\left[\begin{matrix}1, \tfrac{v}{2}+1, \tfrac{v}{2}+\tfrac{3}{2} \\ 2, \tfrac{5}{2}, -\tfrac{v}{2}+\tfrac{3}{2}, \tfrac{v}{2}+\tfrac{1}{2}\end{matrix}; -x^2\right]\right\}$$

and for $e = v+1$, Theorem 3 yields,

$$e^{-x}\sum_{n\geq 0}\frac{(-x)^n}{(\mu)_n}\left\{L_n^{(\nu)}(x)-\tfrac{x}{\nu+2}L_n^{(\nu+1)}(x)\right\}=$$

$$ {}_3F_3\left[\begin{matrix}\nu+3,\tfrac{\mu}{2}+\tfrac{\nu}{2}+1,\tfrac{\mu}{2}+\tfrac{\nu}{2}+\tfrac{1}{2}\\ \nu+2,\mu,\mu+\nu+1\end{matrix};-4x\right]+\frac{2x}{\mu}\,{}_3F_3\left[\begin{matrix}\nu+3,\tfrac{\mu}{2}+\tfrac{\nu}{2}+1,\tfrac{\mu}{2}+\tfrac{\nu}{2}+\tfrac{3}{2}\\ \nu+2,\mu+1,\mu+\nu+2\end{matrix};-4x\right]$$

$$+\frac{x^2}{\mu(\mu+1)}\,{}_3F_3\left[\begin{matrix}\nu+3,\tfrac{\mu}{2}+\tfrac{\nu}{2}+2,\tfrac{\mu}{2}+\tfrac{\nu}{2}+\tfrac{3}{2}\\ \nu+2,\mu+2,\mu+\nu+3\end{matrix};-4x\right]$$

*** Corresponding author: *Yashoverdhan Vyas***

Department of Mathematics, School of Engineering,
Sir Padampat Singhania University, Udaipur, Rajasthan.

*Email: yashoverdhan.vyas@spsu.ac.in*